\documentclass[a4paper,12pt,fleqn]{article}

\usepackage{amssymb}
\usepackage{amsmath}
\usepackage{graphicx}
\usepackage{color}
\usepackage{eurosym}
\usepackage{txfonts}
\usepackage{courier}
\usepackage{pifont}

\newcommand{\rl}{\mathbb{R}}

\newcommand{\set}[1]{\{#1\}}

\newcommand{\cA}{{\cal A}}

\newtheorem{Thm}{Theorem}

\newtheorem{Lm}[Thm]{Lemma}
\newtheorem{Cr}[Thm]{Corollary}
\newtheorem{Alg}[Thm]{Algorithm}

\author{Carrie Rutherford\thanks{London South Bank University}, Robin
Whitty\thanks{Queen Mary University of London}}
\date{\today}
\title{An Algorithmic Proof of the Piff--Welsh Theorem on Transversal Matroid Representations}

\begin{document}

\maketitle

\hspace{.25in}
\parbox[c]{6in}{\small {\bf Abstract} A fundamental theorem of matroid theory establishes that a transversal matroid is representable over fields of any characteristic. It was proved in 1970 by Piff and Welsh: their proof is elegant and concise and, moveover, constructive. However it is far from being algorithmic, in terms of suggesting a step-by-step procedure for deriving a collection of vectors over a given base field representing the transversal matroid induced by a given set system. In this note we recast Piff and Welsh's proof in algorithmic form.
}

{\bf Keywords:} matroid theory, transversal theory, set system, system of distinct representatives.

\section{Background theory}
A {\em set system} consists of $S$, a finite set (the {\em ground set}), together with a finite collection, $\cA$, of subsets of~$S$. A {\em partial transversal} of $\cA$ is a subset of $S$ whose members may be placed in one-to-one correspondence with members of $\cA$ to which they, respectively, belong. A {\em (complete) transversal} (or {\em system of distinct representatives}) is a partial transversal whose cardinality matches that of $\cA$. Transversals are a ubiquitous concept in optimisation: $S$ is a set of contractors and $\cA$ a collection of projects represented by the subsets of contractors eligible for the projects, for example.

For a concrete visualisation, in what follows, we take a small example set system over ground set $S=\{a,b,c,d,e\}$, represented by its {\em incidence matrix}. The rows of the matrix are indexed by the sets in $\cA$ and the columns are indexed by the elements of $S$: a one (zero) in row~$i$, column $j$, indicates that set~$i$ contains (does not contain) element~$j$ of~$S$. In our example~$\cA$ consists of three sets,
$A=\set{a,c,e},\ B=\set{a,b,d,e}$ and $C=\set{a, b,d}$ and the incidence matrix is given as
$$\bordermatrix{& a & b & c & d & e\cr
A & 1 & 0 & 1 & 0 & 1\cr
B & 1 & 1 & 0 & 1 & 1\cr
C & 1 & 1 & 0 & 1 & 0}.
$$
The subset $\{a,b\}$ of $S$ is a partial transversal for this system by virtue of the one-to-one correspondence $a\rightarrow A,\,b\rightarrow B$. This extends to a complete transversal $\{a,b,c\}$ but notice we must `backtrack' to justify this with a one-to-one correspondence: $a\rightarrow B,\,b\rightarrow C,\,c\rightarrow A$.

As our example confirms, the one-to-one correspondence justifying a (partial) transversal cannot always be constructed greedily. However, transversals themselves, i.e. subsets of $S$, can be constructed greedily: if a partial transversal of cardinality $k$ is possible then it is guaranteed that any partial transversal of cardinality $k-1$ can be augmented with an unused element of the ground set $S$. Moreover, if the elements of $S$ have costs attached then this guarantee holds good if we always choose a least-cost unused element of $S$ for augmentation. The augmentation property makes the set of partial transversals of a set system into a {\em matroid}; the cost-insensitivity of augmentation is what makes this matroid important in optimisation. It is interesting that we can identify (partial) subsets of eligible contractors who can be matched to a collection of projects; it is much more interesting that we can do this while always matching the cheapest contractor available for their respective job. (We reiterate that we are only talking about {\em existence} of the cheapest matching! The matching itself, i.e. the one-to-one correspondence justifying the existence claim, although not constructible greedily, can be found in quadratic time via, say, the intersection of matroids).

Thus, associated with a set system $(S, \cA)$ we have a matroid $T(\cA)$ (over ground set $S$) consisting of the set of all partial transversals of $\cA$ over $S$. The partial transversals are called the {\em independent sets} of $T$. Their essential property, corresponding to the guarantee of augmentability, is that `maximal' is equal to `maximum', whereby all maximal (non-augmentable) independent sets  have the same (maximum) cardinality. So all maximal partial transversals have the same size, which may or may not be the cardinality of  $\cA$.

\section{Representations of matroids}
The prototype of the idea of independent sets is the property of linear independence among collections of vectors over a field $K$, say, the real numbers. We say that a matroid $M$ over ground set $S$ is {\em represented} by a collection $C$ of vectors over a field $K$, not necessarily distinct but all having the same length, if the elements of~$S$ can be put into one-to-one correspondence with the vectors of $C$ in such a way that a subset of ground set elements is independent in $M$ if and only if the corresponding vectors are linearly independent over $K$. Normally, we think of the vectors as being the columns of a matrix. Representability of matroids is fundamental to the theory of matroids, supporting classifications and providing powerful proof techniques. But it is also a great asset in practical applications of matroids to be able to work in terms of linear algebra. For example, the  intersection of two matroids (mentioned above in terms of finding matchings) becomes multiplication of the matrices representing them.

It is not obvious that transversal matroids are representable over $\rl$, say. In our example above, the incidence matrix immediately supplies a column vector for each element of $\{a,b,c,d,e\}$. These vectors fail, however, to represent the transversal matroid. The set $\{a,b,c\}$ is a transversal (a maximum independent set) and the first three vectors of the matrix are indeed linearly independent. However, the set $\{a,b,d\}$ is also a transversal (witness the correspondence $a\rightarrow A,\,b\rightarrow B,\,d\rightarrow C$) but columns $a,b,d$ of the incidence matrix fail to be linearly independent. This transversal matroid is, nevertheless, representable over $\rl$: it is sufficient to replace two entries in the incidence matrix with entries other than 1. How exactly we effect this replacement is the subject of this note.

For a special type of set system the incidence matrix is guaranteed to represent the transversal matroid: if the sets of $\cA$ partition $S$ into disjoint subsets then selecting entries for a partial transversal is identical to choosing mutually distinct columns of the incidence matrix. This observation will justify the initialisation step of our algorithm so it is worth stating it formally. In the case where $\cA$ partitions~$S$, the transversal matroid is referred to as a {\bf partition matroid}.

\begin{Lm}
\label{thm:partition}
A partition matroid is represented over any field by its incidence matrix.
\end{Lm}

For the example set system we gave earlier, the initialisation step of our algorithm will extract the following incidence matrix of a partition of $\{a,b,c,d,e\}$:
$$\bordermatrix{& a & b & c & d & e\cr
A & 1 & 0 & 1 & 0 & 1\cr
B & 0 & 1 & 0 & 1 & 0\cr
C & 0 & 0 & 0 & 0 & 0}.
$$
The sets $A$, $B$ and $C$ are different from before but we have retained the names because the idea of the algorithm will be to restore the original sets by reinstating the missing elements progressively. As it is written down here, this incidence matrix represents (over any field) the partition matroid whose maximum independent sets we may list explicitly as: $\{a,b\},\,\{a,d\},\,\{b,c\},\,\{b,e\},\,\{c,d\},\,\{d,e\}$.

Although in general not all transversals correspond to linearly independent sets of columns in an incidence matrix, the reverse does hold: all non-transversals are guaranteed to correspond to linearly dependent sets of columns. This is the content of various equivalent basic theorems in combinatorial optimisation: the marriage theorem, the Frobenius--K\H{o}nig theorem, etc which, informally, say that submatrices of zeros are the only obstacles to the construction of transversals. To be more formal, we introduce another equivalent result which concerns transversal matroid representation. Suppose we take the incidence matrix of a set system but,
instead of a `$+1$' in position $ij$, we place a real number $t_{ij}$. Suppose the totality of the $t_{ij}$ are chosen to form an algebraically independent set. Then we have:
\begin{Thm}[Mirsky and Perfect\cite{Mirsky}]
\label{thm:transrepI}\index{representation!of transversal matroid}
Let $S=\set{a_1,\ldots, a_n}$ be a finite set  and
 $\cA=\set{A_1,\ldots, A_m}$  a collection of subsets of $S$.
Let $X_{\cA}$ be the $m\times n$ matrix whose $ij$-th entry
$x_{ij}$ is defined by
$$x_{ij}=\left\{\begin{array}{lll}t_{ij}&&\mbox{if }A_i \mbox{ contains
}a_j\\0&&\mbox{otherwise}\end{array}\right.,$$ where the $t_{ij}$
are an algebraically independent set. Then $X_{\cA}$ represents the transversal
matroid $T(\cA)$ over $\rl$. In particular, any transversal matroid is representable over the real numbers.
\end{Thm}

We may think of the entries $t_{ij}$ in the Mirsky--Perfect representation of transversal matroids as indeterminates, rather than real numbers. Although mathematically this is a big difference, algorithmically it amounts to the same thing (indeed, it brings the representation within the reach of computer algebra, although not for realistically sized set systems).

\section{The Piff--Welsh Theorem}

\label{subs:piffwelsh} While the Mirsky-Perfect representation is
elegant and revealing it is unsatisfactory to have to work with algebraically independent sets, or with
indeterminates, instead of, say, real numbers or binary numbers, not least because
it puts out of reach much of computational linear algebra.  Luckily,
Dominic Welsh at Oxford University and his DPhil student Mike Piff
simplified matters dramatically in 1970 by proving that the
indeterminates in the Mirsky-Perfect representation may be
replaced by suitably chosen values from any
sufficiently large field:

\begin{Thm}[The Piff--Welsh Theorem\cite{Piff}]
\label{thm:PiffWelsh}
If $M$ is a transversal matroid then $M$ can be represented over  sufficiently large  fields of any characteristic.
\end{Thm}
In particular transversal matroids are representable over infinite (characteristic zero) fields such as the rational numbers. Indeed, since  we can multiply
 vectors in a representation by scalars this means any transversal matroid can be represented by an integer matrix.

Piff and Welsh's original proof of representability (see \cite{Piff} or \cite{Welsh}) relies on a lemma which preserves representability of a matroid under a surjective mapping of its ground set on to a smaller ground set. Essentially the lemma shows that a representation for a larger  matroid, having less `structure', can be reduced to a representation for a smaller matroid with more structure by merging pairs of ground set elements, together with their representing vectors.  We shall refer to this as the Piff--Welsh Merge Algorithm.

In a helpful re-interpretation of the proof given by Oxley\cite{Oxley}, the effect of Piff--Welsh merging is viewed as splitting a given set system into a partition of its ground set augmented with new indeterminate elements $t_{ij}$. As we saw, we can represent a partition matroid directly by its incidence matrix, so this is a larger but very uncomplicated matroid. The main work in constructing the Piff--Welsh representation involves progressively assigning values to the indeterminates $t_{ij}$. Assigning a value to $t_{ij}$ has the effect of merging it with the $i$-th element of the original ground set; it is this merging which inserts the assigned value back into the original incidence matrix.

The main contribution of our presentation is to recast Piff and Welsh's lemma, as viewed by Oxley, as the application of elementary linear algebra to the submatrices of a `tableau'. The proof becomes a systematic enumeration of determinants of submatrices with the complete enumeration then indicating a suitable value for the next `unknown' indeterminate $t_{ij}$. Before stating the algorithm formally we will present it by applying it to our example set system.

 On the left of the `equation' below, the set system is $A=\set{a,c,e},\ B=\set{a,b,d,e}$ and $C=\set{a, b,d}$ and the Mirsky-Perfect representation of its transversal matroid is shown; on the right $a$ has been `split off' from sets $B$ and $C$ as the indeterminates $t_{21}$ and $t_{31}$. Similarly, $b$ has been split off from $C$, $d$ has been split off from $C$ and $e$ has been split off from $B$. The result is a partition of an augmented ground set: $\set{a,b,c,d,e,t_{21},t_{31},t_{32},t_{34},t_{25}}$. The ten columns are a representation of the corresponding partition matroid.
\begin{equation}
\label{eqn:piffwelsh}
\bordermatrix{& a & b & c & d & e\cr
A & t_{11} & 0 & t_{13} & 0 & t_{35}\cr
B & t_{21} & t_{22} & 0 & t_{24} & t_{25}\cr
C & t_{31} & t_{32} & 0 & t_{34} & 0}
\ =\
\bordermatrix{& a & b & c & d & e\cr
A & 1 & 0 & 1 & 0 & 1\cr
B & 0 & 1 & 0 & 1 & 0\cr
C & 0 & 0 & 0 & 0 & 0}
\hspace{.1in}\oplus\hspace{.1in}
\bordermatrix{& t_{21} & t_{31} & t_{32} & t_{34} & t_{25}\cr
A & 0 & 0 & 0 & 0 & 0\cr
B & 1 & 0 & 0 & 0 & 1\cr
C & 0 & 1 & 1 & 1 & 0}
\end{equation}
We shall show how the Piff--Welsh merge algorithm gives new values to the indeterminates $t_{ij}$ by merging back all the indeterminate columns created in the split; some of these new values must clearly be other than +1.
Once we have seen it in action, the algorithm will be stated formally.

The validity of our merge step depends on the familiar fact that the determinant function is $n$-linear in the columns of the matrix. Specifically:
\begin{Lm}
\label{lm:nlinear} Suppose $X$ is an $n\times n$ matrix whose
$i$-th row is the vector $v$ and suppose
that $v$ is expressed as a sum of two vectors: $v=v_1+v_2$. For
$j=1,2$, let $X_j$ denote $X$ with first row replaced with $v_j$.
Then
$$\det X=\det X_1+\det X_2.$$
\end{Lm}

To present the application of this lemma as systematically as possible, we set out the vectors from the right-hand side of equation~(\ref{eqn:piffwelsh}) in the form of a `tableau'. Our first insertion will insert the $t_{21}$ entry  into the $a$ column, so these columns of the tableau are paired and highlighted:
$$
\hspace{-1in}
\begin{array}{r|rr|rrrrrrrr|}\cline{2-3}
\mbox{\bf Iteration 1: indeterminate $t_{21}$} \hspace{.1in} & a & t_{21} & b & c & d & e & t_{31} & t_{32} & t_{34} & \multicolumn{1}{r}{t_{25}}\\\cline{2-11}
&1 & 0 & 0 & 1 & 0 & 1 & 0      & 0      & 0      & 0  \\
&0 & 1 & 1 & 0 & 1 & 0 & 0      & 0      & 0      & 1  \\
&0 & 0 & 0 & 0 & 0 & 0 & 1      & 1      & 1      & 0  \\\cline{2-11}
& & &  & &  &  & & & & \multicolumn{1}{r}{} \\\cline{2-3}
\multicolumn{1}{r}{}& \multicolumn{1}{c}{L} &   \multicolumn{1}{r}{R}  &    &  & &  &      &         &       &\multicolumn{1}{r}{}\\
  \end{array}
$$
Beneath the highlighted columns we are going to construct column vectors called $L$ and $R$; for brevity we shall also use $L$ and $R$ to refer to the highlighted columns (so, above, column $a$ is the `$L$' column and $t_{21}$ is the `$R$' column --- notice we are also abusing notation by using set element names when referring to vectors, thus `column $a$' is short for `the column indexed by $a$')). We will set out a tableau  like the one above for each $t_{ij}$  that has to be inserted into columns $a$--$e$. Sometimes we can just use the value $t_{ij}=1$ but not always: the purpose for the tableau is to identify combinations of columns which {\em prevent us} from making the assignment $t_{ij}=1$. This information will be placed in the $L$ and $R$ vectors below the tableau. What are these combinations of columns? The answer is as follows, being the first of  two key ingredients to getting the final value for each $t_{ij}$:
\begin{center}
\framebox[145mm][c]{\parbox[c]{140mm}{
\begin{description}
\item[Step I:]
 Among the columns other than $L$ and $R$, find a collection $C$ of columns which combines with the $L$ column to give a square non-singular matrix (call it $X_L$) and also combines with the $R$ column to give a square non-singular matrix (call it $X_R$). Record $\det X_L$ in the $L$ vector and  $\det X_R$ in the $R$ vector.
\end{description}
}}
\end{center}
We repeat this for all collections of columns which qualify by giving pairs of non-singular matrices. It is easy to see that, for an initial tableau like the above, it is impossible to satisfy the condition of step~I: because each column has a single nonzero entry, any nonsingular square matrix formed using the L column must become singular when column~L is replaced by column~R.

If no pairs of non-singular matrices can be formed then there is nothing to prevent us from setting $t_{ij}=1$. So for our first insertion we have found we can place $t_{21}$ back into column $a$ of our matrix with value 1. This gives us a new tableau in which we are seeking to merge indeterminate $t_{31}$ into column $a$:

$$
\hspace{-1in}
\begin{array}{r|rr|rrrrrrr|}\cline{2-3}
\mbox{\bf Iteration 2: indeterminate $t_{31}$} \hspace{.1in} & a & t_{31} & b & c & d & e & t_{32} & t_{34} & \multicolumn{1}{r}{t_{25}}\\\cline{2-10}
&1 & 0 & 0 & 1 & 0 & 1 &  0      & 0      & 0  \\
&1 & 0 & 1 & 0 & 1 & 0 &  0      & 0      & 1  \\
&0 & 1 & 0 & 0 & 0 & 0 &  1      & 1      & 0  \\\cline{2-10}
& & &  & &     & & & & \multicolumn{1}{r}{} \\\cline{2-3}
\multicolumn{1}{r}{}& \multicolumn{1}{c}{L} &   \multicolumn{1}{r}{R}  &    &   &  &      &         &       &\multicolumn{1}{r}{}\\
  \end{array}
$$
Notice that column $t_{21}$ and column $a$ have been merged in this new tableau. This means that we  no longer have a partition matroid. However, the incidence matrix for the set system over the reduced ground set of nine elements {\em is still a representation for the associated transversal matroid}. This is because, in the current case vacuously but generally by invoking  Lemma~\ref{lm:nlinear}, any nonzero determinants for matrices involving columns L and R will sum, in the merged matrix, to again give a nonzero determinant. Meanwhile, we cannot create new nonzero determinants by virtue of the Mirsky--Perfect theorem.

Before repeating our search for nonsingular matrices we identify a short-cut: it is unnecessary to include in our tableau any repeated columns since such columns can never appear together in a nonsingular matrix. So the above tableau becomes
$$
\hspace{-1in}
\begin{array}{r|rr|rr|}\cline{2-3}
\mbox{\bf Iteration 2: indeterminate $t_{31}$} \hspace{.1in} & a & t_{31} & b & \multicolumn{1}{r}{c}\\\cline{2-5}
&1 & 0 & 0 & 1    \\
&1 & 0 & 1 & 0    \\
&0 & 1 & 0 & 0   \\\cline{2-5}
& & &  &   \multicolumn{1}{r}{} \\\cline{2-3}
\multicolumn{1}{r}{}& \multicolumn{1}{c}{L} &   \multicolumn{1}{r}{R}  &           &\multicolumn{1}{r}{}\\
  \end{array}
$$

There is again no collection of two columns in the new tableau which we can combine with both the $L$ column (i.e. $a$) and the $R$ column (i.e. $t_{31}$) to get non-singular matrices. This means there is again nothing to prevent us inserting $t_{31}$ into $a$ using value 1. So now column $t_{31}$ is surpressed and the final entry in column $a$ becomes 1.

We are ready to make our third insertion, combining $t_{32}$ into the $b$ column. Omitting repeated columns we have:

$$
\hspace{-1in}
\begin{array}{r|r|rr|r|}\cline{3-4}
\multicolumn{1}{r}{\mbox{\bf Iteration 3: indeterminate $t_{32}$} \hspace{.1in} }&  \multicolumn{1}{r|}{a} & b & t_{32} & \multicolumn{1}{r}{c}\\\cline{2-5}
&1 & 0 & 0 & 1    \\
&1 & 1 & 0 & 0    \\
&1 & 0 & 1 & 0   \\\cline{2-5}
\multicolumn{1}{r}{}& \times  & -1 &   1  &       \multicolumn{1}{r}{\times}\\\cline{3-4}
\multicolumn{2}{r}{}  & \multicolumn{1}{c}{L} &   \multicolumn{1}{r}{R}  &          \multicolumn{1}{r}{}
  \end{array}
$$
The submatrix of columns $a,b,c$ has determinant $-1$, recorded in the $L$ column; the submatrix of columns $a,t_{32},c$ has determinant $1$, recorded in the $R$ column. The columns involved in the nonsingular matrices in each row of the $L$ and $R$ vectors, in this case the unique choice of $a$ and $c$, are indicated by $\times$ signs.

Now observe that, in invoking Lemma~\ref{lm:nlinear}, we cannot simply place a 1 in the last entry of column~$b$ because the sum of the two determinants would be zero. At this point we introduce the
  second key element in our algorithm:
\begin{center}
\framebox[145mm][c]{\parbox[c]{140mm}{
\begin{description}
\item[Step II:] Assign to $t_{ij}$ any non-zero value such that the vector sum $L+t_{ij}\times R$ has no zero entries.
    \end{description}
    }}
\end{center}

For the above tableau, we take our single-entry vectors $L=(-1)$ and $R=(1)$ and observe that $t_{23}=-1$ will give
$$L+t_{23}\times R=(-1)+\,-1\times(1)=(-1)+(-1)=(-2),$$
a vector whose single entry is non-zero, as required. Note that $t_{23}=0$ would also have given a non-zero sum but the $t_{ij}$ must be chosen to be non-zero, to comply with the Mirsky--Perfect theorem (Theorem~\ref{thm:transrepI}). However, $t_{23}=2$ would have been another valid choice.

We now come to our fourth insertion, as we place $t_{34}$ into the $d$ column of our representation. Notice that the $b$ column has  had the valid value $t_{23}=-1$ inserted in the final entry:
$$
\hspace{-1in}
\begin{array}{r|rrr|rr|}\cline{5-6}
\multicolumn{1}{r}{\mbox{\bf Iteration 4: indeterminate $t_{34}$} \hspace{.1in} }&  \multicolumn{1}{r}{a} & b & c & d & t_{34}\\\cline{2-6}
&1 & 0  & 1 & 0  &0  \\
&1 & 1  & 0 & 1  &0   \\
&1 & -1 & 0 & 0  &1 \\\cline{2-6}
\multicolumn{1}{r}{}& \times & \times & &   1  & 1 \\
\multicolumn{1}{r}{}& \times &  & \times &  1  & -1 \\
\multicolumn{1}{r}{}&  & \times & \times &  -1  & -1 \\\cline{5-6}
\multicolumn{4}{r}{}  & \multicolumn{1}{c}{L} &   \multicolumn{1}{r}{R}
  \end{array}
$$

This time we may make the assignment $t_{34}=2$ to ensure $L+t_{34}\times R$ has no nonzero entries, so that Lemma~\ref{lm:nlinear} applies correctly to the merging of columns $d$ and $t_{34}$.

In or final tableau, combining columns $e$ and $t_{25}$, we omit column $c$ which is a repeat of column $e$:
$$
\hspace{-1in}
\begin{array}{r|rrr|rr|}\cline{5-6}
\multicolumn{1}{r}{\mbox{\bf Iteration 5: indeterminate $t_{25}$} \hspace{.1in} }&  \multicolumn{1}{r}{a} & b & d & e & t_{25}\\\cline{2-6}
&1 & 0  &   0  &1 &0 \\
&1 & 1  &   1  &0 & 1  \\
&1 & -1 &   2  &0 & 0 \\\cline{2-6}
\multicolumn{1}{r}{}& \times & \times & &  -2  & 1 \\
\multicolumn{1}{r}{}& \times &  & \times &  1  & -2 \\\cline{5-6}
\multicolumn{4}{r}{}  & \multicolumn{1}{c}{L} &   \multicolumn{1}{r}{R}
  \end{array}
$$
A value of $t_{25}=1$ is valid here. So we have our  completed representation for the transversal matroid of the set system $A=\set{a,c,e},\ B=\set{a,b,c,d}$ and $C=\set{a,b,d}$:
$$\bordermatrix{& a & b & c & d & e\cr
A & 1 & 0 & 1 & 0 & 1\cr
B & 1 & 1 & 0 & 1 & 1\cr
C & 1 & -1 & 0 & 2 & 0}.
$$

We now give a formal statement of the Piff--Welsh merge algorithm.

\begin{Alg}[The Piff--Welsh Merge Algorithm]\index{Piff--Welsh
merge}\index{algorithm!Piff--Welsh merge}
\label{alg:PiffWelsh}
 \ \vspace{-.1in} \em
 {\em
\begin{tabbing}
\hspace{.25in}\=\hspace{.25in}\=\hspace{.25in}\=\hspace{.25in}\=\hspace{.25in}\=\hspace{.25in}\=\kill
\>{\bf INPUT:}\>\>\> $m\times n$ matrix $X$\\
\>\>\>\>vectors $v_L$ and $v_R$ chosen from the columns of $X$\\
\>{\bf OUTPUT:}\>\>\>real number $t$\\
{\em 1}\>\>${\cal C}$:=set of all pairs $(X_L, X_R)$ of $m\times m$ non-singular matrices of the form $X_L=[v_L\,|\,Y]$\\
 \>\>\> and $X_R=[v_R\,|\,Y]$, where $Y$ is an $m\times (m-1)$ submatrix of $X$ not including columns $v_L$ or $v_R$;\\
{\em 2}\>\>$Z$:=a matrix with two cols, initially having zero rows;\\
{\em 3}\>\>{\bf for} each pair $(X_L, X_R)\in {\cal C}$ {\bf do}\\
{\em 4}\>\>\>append the row vector $(\det X_L,\det X_R)$ to $Z$;\\
{\em 5}\>\>{\bf od};\\
{\em 6}\>\>{\bf return} a real number $t$ such
that $Z\times (1,t)^T$
has no zero elements
\end{tabbing}
}
\end{Alg}

We have argued the correctness of this algorithm informally during the presentation of our running example above. To give a formal argument for the merge algorithm, observe firstly that step~6 in
Algorithm~\ref{alg:PiffWelsh} is always possible, provided we have
sufficient choices for~$t$. This is certainly the case
for $\rl$; in the case of a finite field of a given prime characteristic $p$ it is sufficient to work over the field $GF(p^k)$ for high enough $k$. Secondly we verify that the choice of $t$ in step~6 ensures that the algorithm preserves independent sets.
Suppose that the vectors in $X$ are a representation and that we want to merge $v_L$
and $v_R$ to represent the same ground set element, and suppose $X_L$ corresponds to a set of linearly independent columns (so it is a non-singular submatrix). Replace $v_L$ in
$X_L$ with $v_L+tv_R$ to get $X'_L$. Then by
Lemma~\ref{lm:nlinear}, $\det X'_L=\det X_L+t\det
X_R=(\det X_L,\det X_R)\times(1,t)^T$, which is non-zero by the choice of $t$ in step~6. So $X'_L$ is again nonsingular. Conversely, dependent subsets of $X$ are also preserved
by Algorithm~\ref{alg:PiffWelsh} due to the the Mirsky--Perfect theorem (Theorem~\ref{thm:transrepI}).

The merge algorithm applies to a matrix $X$. In our proof of the algorithm, $X$ is the whole tableau, comprising all original ground set elements and invariants $t_{ij}$. In our running example we reduced the tableaux and, in fact, no invariants beyond the one currently being merged were included in any tableaux. It might be thought that the algorithm need not `look ahead' to future invariants but this is not in general the case: columns of $X$ may be omitted only where repeating columns would produce singular submatrices (step~1 in the algorithm). Consider, for example, the following partition matroid on five elements, of which three are invariants:
$$\bordermatrix{& a & b & t_{21} & t_{22} & t_{32}\cr
A & 1 & 1 & 0 & 0 & 0\cr
B & 0 & 0 & 1 & 1 & 0\cr
C & 0 & 0 & 0 & 0 & 1}.
$$
After $t_{21}$ is merged into the `$a$' column, automatically taking the value $t_{21}=1$, it is the turn of $t_{22}$ to be merged into the `$b$' column. In the corresponding tableau we must include the $t_{32}$ column, leading to the evaluation $t_{22}\neq 1$. Without $t_{32}$ there would trivially be no nonsingular matrices, giving an incorrect evaluation $t_{22}=1$. Of course, in the final representation of the transversal matroid on ground set $\{a,b\}$ there are in any case no $3\times 3$ nonsingular submatrices! But a value of $t_{22}=1$ would imply incorrectly that $\{a,b\}$ is not a transversal for sets $A$ and $B$.

We conclude by mentioning the fact that, as well as partition matroids, another class of matroids is representable thanks to the Piff--Welsh theorem: recall the uniform matroid $U_{k,n}$ is the matroid in which the independent sets are all subsets of size at most $k$ of a ground set of size $n$.
\begin{Cr}
\label{cr:PiffWelsh}
Any uniform matroid may be represented by a matrix with integer entries.
\end{Cr}

\medskip
{\bf Proof.} Suppose the ground set of $U_{k,n}$ is $A=\set{a_1, \ldots, a_n}$. Let $\cA$ be the set system which consists of $k$ copies of $A$. Then any subset of size $k$ of $A$ is a transversal for $\cA$ and vice versa. So $T(\cA)$ is isomorphic\index{matroid!isomorphism} to $U_{k,n}$. Since by Theorem~\ref{thm:PiffWelsh} $T(\cA)$ is representable over the field of rational numbers, we may represent $U_{k,n}$ by a matrix whose entries are rationals. It suffices to multiply this matrix by the least common multiple of the denominators of its entries to get a representation using only integers.

\end{document}